\documentclass{amsart}
\usepackage{amsmath}%
\usepackage{amsfonts}%
\usepackage{amssymb}%
\usepackage{graphicx}
\theoremstyle{plain}

\numberwithin{equation}{section}
\begin{document}
\title[ variation of dynamical canonical heights, Intersection numbers]{ On  variation of dynamical canonical heights, and Intersection numbers}
\author{Jorge Mello}
\address[]
{Universidade Federal do Rio de Janeiro, Instituto de Matem\'{a}tica. mailing adress:\newline Rua Aurora 57/101, Penha Circular, 21020-380 Rio de Janeiro, RJ, Brasil. %
 } \email[]{jmelloguitar@gmail.com}%
\urladdr{https://sites.google.com/site/algebraufrj/}

\thanks{}
\date{July 28, 2017}
\subjclass{} %
\keywords{Canonical Heights, Local heights, Intersection numbers, Families of varieties, Dynamical systems.}%
\dedicatory{}

\begin{abstract}
 We study families of varieties endowed with polarized canonical eigensystems of several maps, inducing canonical heights on the dominating variety as well as on the "good" fibers of the family. We show explicitely the dependence on the parameter for global and local canonical heights defined by Kawaguchi when the fibers change, extending previous works of J. Silverman and others. Finally, fixing an absolute value $v \in K$ and a variety $V/K$, we descript the Kawaguchi`s canonical local height $\hat{\lambda}_{V,E,\mathcal{Q},}(.,v)$ as an intersection number, provided that the polarized system  $(V,\mathcal{Q})$ has a certain weak N\'{e}ron model over Spec$(\mathcal{O}_v)$ to be defined and under some conditions depending on the special fiber.  With this we extend N\'{e}ron's work strengthening Silverman's results, which were for systems having only one map.
\end{abstract}
\maketitle

\section{Introduction}
Over abelian varieties $A$ defined over a number field $K$, N\'{e}ron and Tate constructed canonical height functions $\hat{h}_L: A(\bar{K}) \rightarrow \mathbb{R}$ with respect to symmetric ample line bundles $L$ which enjoy nice properties, and can be used to prove Mordell Weil theorem for the rational points of the variety. More generally, in [9], Call and Silverman constructed canonical height functions on projective varieties $X$ defined over a number field which admit a morphism $f:X \rightarrow X$ with $f^*(L) \cong L^{\otimes d}$ for some line bundle $L$ and some $d >1$.
In another direction, Silverman [31] constructed canonical height functions on certain $K3$ surfaces $S$ with two involutions $\sigma_1, \sigma_2$ (called Wheler's $K3$ surfaces) and developed an arithmetic theory analogous to the arithmetic theory on abelian varieties. In another big step, S. Zhang [37] considered a projective variety with a metrized line bundle, and then could define  and evaluate a new canonical height on any subvariety of the variety, not only on points. He developed this using arithmetic intersection theory arising from Arakelov Geometry.

It was an idea of Kawaguchi [16] to consider polarized dynamical systems of several maps, namely, given $X/K$ a projective variety, $f_1,...f_k:X \rightarrow X$ morphisms on defined over $K$, $\mathcal{L}$ an invertible sheaf on $X$ and a real number $d>k$ so that $f_1^*\mathcal{L}\otimes ... \otimes f_k^*\mathcal{L} \cong \mathcal{L}^{\otimes d}$, he constructed a canonical height function associated to the polarized dynamical system $(X, f_1,..., f_k, \mathcal{L})$ that generalizes the earlier constructions mentioned above. In the Wheler's K3 surfaces' case above, for example, the canonical height defined by Silverman arises from the system formed by $(\sigma_1, \sigma_2)$ by Kawaguchi's method. He also used Arakelov geometry to calculate his canonical height of subvarieties, and not only on points, generalizing Zhang's construction.

As a canonical height, it is known that Kawaguchi's height is, up to a constant, equal to a Weil height. We start studying how such bounds can vary more explicitely in families of varieties, when for each fiber we associate one canonical height. This kind of research was made by Silverman and Tate in [30] for families of abelian varieties, and afterwards by Silverman in [9] for a family of general varieties and the canonical height developed there. We therefore generalize their results for the Kawaguchi situation. Namely, for a family $\pi: \mathcal{V} \rightarrow T$ of varieties with a system $\mathcal{Q}$ of maps $\phi_i:\mathcal{V}/T \dashrightarrow \mathcal{V}/T$ and a divisor $\eta$ satisfying $\pi(|(\bigotimes_{i=1}^k \phi_i^* \eta) - \alpha \eta|) \neq T$. Then on all fibers $\mathcal{V}_t$ for $t$ in a certain $T^0$, there is a canonical height $\hat{h}_{\mathcal{V}_t, \eta_t, \mathcal{Q}_t}$, and we ask to bound the difference between this height and a given Weil height $h_{\mathcal{V}, \eta}$ in terms of the parameter $t$. We show that there exist constants $c_1$ and $c_2$ such that
\begin{center}
$|\hat{h}_{\mathcal{V}_t, \eta_t, (\mathcal{Q})_t}(x)- h_{\mathcal{V}, \eta}(x)| \leq c_1h_T(t)+ c_2$~ for all $t \in T^0$ and all $x \in \mathcal{V}_t.$
\end{center}
 We can find also a kind of local version for the above result. In fact, in theorem 4.2.1 of [16] Kawaguchi showed that his canonical height can be seen as a sum of local canonical heights as in the case of abelian varieties. Therefore, we show an estimate for the difference between the canonical local height and a given Weil height.

When $T$ is a curve, $h_T$ is a Weil height associated to a divisor of degree one, $P:T \rightarrow \mathcal{V}$ is a section, $V$ the generic fiber of $\mathcal{V}$ is a variety over the global function field $V(\bar{K}(T))$, $P_t:=P(t)$, and the section $P$ corresponds to a rational point $P_V \in V(\bar{K}(T))$, we show that 
\begin{center}
$\lim_{h_T(t) \rightarrow \infty, t \in T^0(\bar{K})} \dfrac{\hat{h}_{\mathcal{V}_t, \eta_t, (\mathcal{Q})_t}(P_t)}{h_T(t)} =\hat{h}_{V, \eta_V, (\mathcal{Q})_V}(P_V),$
\end{center}
generalizing a result of Silverman in [9].

The work is also place for an analysis of canonical local heights for non-archimedean places in the context of intersection theory. Inspired in the final section of [9], this theme was already discussed for abelian varieties. We work out with system of several maps again. We prove more generally that if $V$ has a model $\mathcal{V}$ over a complete ring $\mathcal{O}_v$ such that every rational point extends to a section and such that $k$ morphisms $\phi_i:V \rightarrow V$ extend to finite morphisms $\Phi_i: \mathcal{V} \rightarrow \mathcal{V}$, then, under some condition on the special fiber, the canonical local height of Kawaguchi is given by an intersection multiplicity on $\mathcal{V}$. For such generalization, we use Frobenius-Perron theory on eigenvalues of matrices.

In the last section, we point out that the admissible metric, defined by Kawaguchi in [16] for a dynamical system of maps associated to a bundle, does not change if we start with another dynamical system that has maps commuting with the maps of the first system, and that are associated with the same divisor by a similar divisorial relation. The canonical measures risen by both systems would then be the same as well, and the main results of [26] are now in a more general setting.

\section{Two variation theorems}
The following notation will be used for this section and the next section.
\begin{itemize}
\item $ K:$ A global field with characteristic 0 and a complete set of proper absolute values satisfying a product formula. We will call such a field a global height field.

\item $M:= M_{\bar{K}}:$ The set of absolute values on $\bar{K}$ extending those on $K$.

\item $ T/K:$  a smooth projective variety.

\item $ h_{T}:$ A fixed Weil height function  on $ T$ associated to an ample divisor, chosen to satisfy $h_T \geq 0.$

\item $ \mathcal{V}/K:$  a smooth projective variety.

\item $\pi :$ a morphism $\pi: \mathcal{V} \rightarrow T$ defined over $K$ whose generic fiber is smooth and geometrically irreducible.

\item $ \phi_i:$  rational maps $\phi_i:\mathcal{V}/T \dashrightarrow \mathcal{V}/T$ defined over $K$ for $i=1,...,k$, such that $\phi_i$ is morphism on the generic fiber of $\mathcal{V}/T.$ Our assumption that $\phi_i$ is on $\mathcal{V}/T$ means that $\pi \circ \phi_i= \pi$.

\item $ \eta:$ A divisor class $\eta \in $ Pic$(\mathcal{V}) \otimes \mathbb{R}$ satisfying $\bigotimes_{i=1}^k \phi_i^* \eta = \alpha \eta$ for some real $\alpha >k$.

\item $T^0:$ the subset of $T$ having good fibers in the sense that
\begin{center}
$T^0=\{ t \in T : \mathcal{V}_t ~\mbox{is smooth and} ~ (\phi_i)_t: \mathcal{V}_t \rightarrow \mathcal{V}_t ~ \mbox{is a morphism} ~ \forall i \}.$
\end{center}
where $\mathcal{V}_t := \pi^{-1}(t).$

\item  $ \mathcal{Q}_n$ for $ n \in \mathbb{N} :$ the sets of iterates of functions defined as $\mathcal{Q}_0=\{ $Id$\}, \newline \mathcal{Q}_1= \mathcal{Q} =\{\phi_1,...,\phi_k\}$, and $\mathcal{Q}_n=\{\phi_{i_1} \circ ... \circ \phi_{i_n} ; i_j =1,...,k \}.$

\item  $ (\mathcal{Q}_n)_t$ for $ n \in \mathbb{N} :$ the sets of iterates of functions defined as \newline $(\mathcal{Q}_0)_t=\{ $Id$\}, (\mathcal{Q}_1)_t= (\mathcal{Q})_t =\{(\phi_1)_t,...,(\phi_k)_t\}$, and  \newline $(\mathcal{Q}_n)_t=\{(\phi_{i_1})_t \circ ... \circ (\phi_{i_n})_t ; i_j =1,...,k \},$ the restrictions \newline of the $\phi_i'$s to the fiber $\mathcal{V}_t$.
\end{itemize}

We also assume that the divisor class $\bigotimes_{i=1}^k \phi_i^* \eta - \alpha \eta $ is fibral, which means that it can be represented by a divisor $\Delta$ such that $\pi(|\Delta|) \neq T,$ or equivalently, there exists a divisor $D$ in  Div$(T)$ such that $\pi^*D > \Delta >-\pi^*D.$

For any $t \in T^0$ we let $i_t: \mathcal{V}_t \rightarrow \mathcal{V}$ be the natural inclusion, and then by definition $i_t^* \eta=\eta_t.$ The fiber $\mathcal{V}_t$ is irreducible for each $t \in T^0$. If the support of a fibral divisor includes an irreducible fiber, it is always possible to find a linearly equivalent divisor which does not include that fiber. This implies that
\begin{center}
$\bigotimes_{i=1}^k (\phi_i)_t^* \eta_t = \alpha \eta_t \in ~$Pic$(\mathcal{V}_t) \otimes \mathbb{R}$ for all $t \in T^0.$ 
\end{center}
From this and from theorem 1.2.1 of [16], we have that for each $t \in T^0$ there is a canonical height
\begin{center}
$\hat{h}_{\mathcal{V}_t, \eta_t, (\mathcal{Q})_t}: \mathcal{V}_t(\bar{K}) \rightarrow \mathbb{R}.$
\end{center}
We now fix a Weil height 
\begin{center} 
$h_{\mathcal{V}, \eta}: \mathcal{V}(\bar{K}) \rightarrow \mathbb{R}$ 
\end{center} associated to $\eta.$ It follows from the properties of the height functions of Kawaguchi [16], and functoriality of the Weil height function, that
\begin{center}
$\hat{h}_{\mathcal{V}_t, \eta_t, (\mathcal{Q})_t} = h_{\mathcal{V}_t, \eta_t} + O(1)= h_{\mathcal{V}, \eta} \circ i_t + O(1).$
\end{center}
where the $O(1)$ depend on $t$. For $t \in T$, any two canonical heights $\hat{h}_{\mathcal{V}_t, \eta_t, (\mathcal{Q})_t}$ differ from the Weil Height $h_{\mathcal{V}, \eta}$ by a bounded amount constant depending on $t$. For applications, it is important to have an explicit bound for such constant.
Let us check how we can see this dependence explicitly in the following theorem, which is an extension of theorem 3.1 of [9], for one map's systems. Such case was a more general form of the work of Silverman and Tate for families of abelian varieties done in [30]. \newpage

{\bf Theorem 2.1:}
\textit{With notation as above, there exist constants $c_1, c_2$ depending on the family $\mathcal{V} \rightarrow T$, the system $ \mathcal{Q}$, the divisor class $\eta$, and the choice of Weil height functions $ h_{\mathcal{V}, \eta}$ and $h_T$, so that}
\begin{center}
\textit{$|\hat{h}_{\mathcal{V}_t, \eta_t, (\mathcal{Q})_t}(x)- h_{\mathcal{V}, \eta}(x)| \leq c_1h_T(t)+ c_2$~ for all $t \in T^0$ and all $x \in \mathcal{V}_t.$}
\end{center}

\begin{proof}
From the definition of $T^0$, one can conclude that  $\mathcal{V}^0:=\pi^{-1}(T^0)$ is smooth, so we are able to apply  some resolution of singularities to $\mathcal{V}$ without changing $\mathcal{V}^0$. Moreover, although the maps $\phi_i: \mathcal{V} \dashrightarrow \mathcal{V}$ are merely rational, they are morphisms on $\mathcal{V}^0.$ This means we can blow-up $\mathcal{V}$ to produce:

(i) smooth projective varieties $\tilde{\mathcal{V}}_i,$

(ii) birational morphisms $\psi_i: \tilde{\mathcal{V}}_i \rightarrow \mathcal{V}$ which are isomorphisms on $\mathcal{V}^0$

 (iii) morphisms $\xi_i: \tilde{\mathcal{V}}_i \rightarrow \mathcal{V}$ which extend the rational maps  \newline $\phi_i \circ\psi_i: \tilde{\mathcal{V}}_i \rightarrow \mathcal{V}.$ \newline
The existence of $\tilde{\mathcal{V}}_i$ with these properties follows from [13] II.7.17.3 and II.7.16, except that $\tilde{\mathcal{V}}_i$ might be singular. We then use Hironaka's resolution of singularities to make $\tilde{\mathcal{V}}_i$ smooth and we have the desired properties. 

Next we choose a divisor $E \in$ Div$(\mathcal{V}) \otimes \mathbb{R}$ in the divisor class of $\eta$, and we let $H \in $ Div$(T)$ be the ample divisor used to define $h_T.$ Our assumption that $\bigotimes_{i=1}^k \phi_i^* \eta - \alpha \eta$ is fibral guarantees the existence of a divisor $D \in$ Div$(T) \otimes \mathbb{R}$ with
\begin{center}
$\pi^*D > \sum_{i=1}^k \phi_i^* E - \alpha E > -\pi^*D.$
\end{center} 
We also choose an integer $n >0$ so that the divisors
\begin{center}
$nH+D$ and $nH -D$ are both ample on $T$.
\end{center}
The height function with respect to a positive divisor is bounded below out of the support of the divisor, and for an ample divisor such height is everywhere bounded below. So the last assertions imply that
\begin{center}
$|h_{\mathcal{V},\bigotimes_{i=1}^k \phi_i^* E- \alpha E }| \leq nh_{\mathcal{V},\pi^*H} + O(1) = nh_{T,H}\circ \pi + O(1)$  \newline for all points $P$ in $\mathcal{V}^0-|D|.$
\end{center}
 Now let $x \in \mathcal{V}^0$ be any point, and let $\tilde{x}_i \in \tilde{\mathcal{V}_i}$ satisfying $\psi_i(\tilde{x}_i)=x.$  In the following computation, we write $O(1)$ for a quantity that is boundable in terms of the family $\mathcal{V} \rightarrow T$, the maps $\phi_i$, the divisor class $\eta$, and the choice of Weil height functions $h_{\mathcal{V},\eta}=h_{\mathcal{V},E}$ and $h_T=h_{T,H}.$The most important here is that $O(1)$ is independent of $x \in \mathcal{V}^0.$
\begin{center}
$|\sum_ih_{\mathcal{V},E}(\phi_i(x))- \alpha h_{\mathcal{V},E}(x)|~~~~~~~~~~~\:\:\:\:\:\:\:\:\:~~~~~\:\:\:\:\:\:\:\:\:\:\:\:\:\:\:~~~~~\newline \newline
=|\sum_ih_{\mathcal{V},E}((\phi_i\circ \psi_i) (\tilde{x}_i))- \alpha h_{\mathcal{V},E}(x)|  ~~\:\:\:\:\:~~~\newline \newline
=|\sum_ih_{\mathcal{V},E}(\xi_i (\tilde{x}_i))- \alpha h_{\mathcal{V},E}(x)|  ~~~~~~~\:\:\:\:\:\:\:\:\:\:\:\:\:\:\:\:\:~~~~\:~~~~~~\newline \newline
=|\sum_ih_{\mathcal{V},\xi_i^*E}(\tilde{x}_i)- \alpha h_{\mathcal{V},E}(x)| + O(1)  ~~~~~~\:\:\:\:\:\:\:\:~\newline \newline
=|\sum_ih_{\mathcal{V},\psi_i^*\phi_i^*E}(\tilde{x}_i)- \alpha h_{\mathcal{V},E}(x)| + O(1) ~~~\:\:\:\: \newline \newline
=|\sum_ih_{\mathcal{V},\phi_i^*E}(\psi_i(\tilde{x}_i))- \alpha h_{\mathcal{V},E}(x)| + O(1) ~ \newline \newline
=|\sum_ih_{\mathcal{V},\phi_i^*E}(x)- \alpha h_{\mathcal{V},E}(x)| + O(1)  ~~~\:\:\:\:\:\:\:\:~\:~~~~~\newline \newline
=|h_{\mathcal{V},\bigotimes_{i=1}^k \phi_i^* E- \alpha E }| + O(1) ~~\:\:\:~~~~~~~\:\:\:\:\:\:\:\:~\:\:\:\:\:\:\:\:\:\:\:\:\:\:\:\:\:~~~~~~~~~~~~~~~~~\newline \newline
\leq nh_{T,H}(\pi (x)) + O(1).~~~~~~~~~~\:\:\:\:\:\:\:\:\:\:\:\:\:\:\:\:\:\:\:\:\:\:\:\:\:~~~~~\:\:\:\:\:\:\:\:\:\:~~~~~~~~~~~~\:\:\:\:\:\:\:\:\:\:\:\:\:\:\:\:\:\:\:\:\:\:\:\:~~~~~~~~\:~~~~~~~~~~~~~~~~~$
\end{center}
This inequality is valid on $\mathcal{V}^0$ out of the suport of $D$. Choosing different divisors $E$ in the class of $\eta$ to move $D$, we then obtain the inequality for all points in $\mathcal{V}^0$, in a similar way of [30], pages 203-204. This proves

(*) $|\sum_ih_{\mathcal{V},E}(\phi_i(x))- \alpha h_{\mathcal{V},E}(x)| \leq  nh_{T,H}(\pi (x)) + O(1)$ for all $x \in \mathcal{V}^0.$

In order to complete the proof, we remember theorem 1.2.1 of [16], which says in similar conditions that \begin{center} (1): ~ $|\sum_{i=1}^kh_L(f_i(x))-dh_L(x)|\leq C$  implies (2):~ $|\hat{h}_{L,\mathcal{F}}(x)-h_L(x)| \leq \frac{C}{d-k}$. \end{center}We use (*) in place (1) and obtain 
\begin{center}
$|\hat{h}_{\mathcal{V}_t, \eta_t, (\mathcal{Q})_t}(x)- h_{\mathcal{V}, \eta}(x)| \leq \dfrac{ nh_{T,H}(\pi (x)) + O(1)}{\alpha-k}.$
\end{center}
\end{proof}
 By $\lambda_{\mathcal{V},E}$ we denote a Weil local height function \begin{center}$\lambda_{\mathcal{V},E}:(\mathcal{V} - |E|)\times M_{\bar{K}} \rightarrow \mathbb{R}$ \end{center}associated to the divisor $E$. We can now give a similar estimate for the difference between the canonical local height $\hat{\lambda}_{\mathcal{V}_t,E_t,(\mathcal{Q})_t}$ defined by Kawaguchi in theorem 4.2.1 of [16] and a given Weil local height $\lambda_{\mathcal{V},E}$, generalizing Lang`s result [20] for abelian varieties. Here we make extensive use of the notation in theorem 7.3 and corollary 7.4 of [27], and obtain a local version of theorem 2.1. So 2.1 can be again deduced adding the following theorem over all absolute values of $K$ and applying theorem 4.3.1 of [16]. Moreover,  we note that one may skip to use resolution of singularities to prove the previous theorem just proving the next and adding up local contributions, even obtaining a stronger result, since avoiding to use resolution of singularities gives us that theorem 2.1 is true in any characteristic, using the machinery of [27]. So theorem 2.1 is valid any global field $K$. For basic facts about local height functions, $M_K-$bounded functions and $M_K-$ constants, see [20], chapter 10. We here use freely terminology from [27].

{\bf Theorem 2.2:} \textit{With notation as above, fix a divisor} $E \in$ Div$(\mathcal{V}) \otimes \mathbb{R}$ \textit{in the class of $\eta$ and a Weil local height function} $\lambda_{\mathcal{V},E}.$ \textit{Let $U$ be defined as the following set}
  $\{ t \in T: \mathcal{V}_t ~ \mbox{\textit{is smooth}},  (\phi_i)_t: \mathcal{V}_t \rightarrow \mathcal{V}_t \mbox{ \textit{are morphisms} }, E_t \mbox{ \textit{is a divisor on}}  ~\mathcal{V}_t, \newline ~~~~~~~~~~~~~~~~~~~~~~~~~~~~~~~~~~~~~~~~~~~~~~~~~~~~~~~~~~~~~~~~\mbox{ \textit{and} } \sum_i(\phi_i)_t^*E_t \sim \alpha E_t\}.$ \newline \newline
\textit{(The condition that $E_t$ be a divisor on $\mathcal{V}_t$ means that $|E|$ contains no component of $\mathcal{V}_t$, see [13], III.9.8.5.)}

\textit{Let $\partial U := T -U$ be the complement of $U$, and let $\lambda_{\partial U}$ be a local height function associated to $\partial U$ as described in [27].}

\textit{It is possible to choose canonical local heights $\hat{\lambda}_{\mathcal{V}_t, E_t,(\mathcal{Q})_t }$ as described in Theorem 4.2.1 of [16], one for each $t \in U$, in such a way that} \newline \newline
$|\hat{\lambda}_{\mathcal{V}_t, E_t,(\mathcal{Q})_t }(x,v) - \lambda_{\mathcal{V},E}(x,v)| \leq c \lambda_{\partial U}(t,v) \newline ~~~~~~~~~~~~~~~~~~~~~~~~~~~~~~~~~$ \textit{for all $(x,v) \in (\mathcal{V}-|E|) \times M$ with $\pi(x)=t \in U.$}

\begin{proof}
We substitute $\mathcal{V}$ by the quasi projective variety $\pi^{-1}(U)$, and substitute $E$ by its restriction to this new $\mathcal{V}.$ This does not affect the statement of the theorem because [27] section 5 says that our old $\lambda_{\mathcal{V},E}$ and our new $\lambda_{\mathcal{V},E}$ differ by $O(\lambda_{\partial U})$. From the definition of $U$ we have that $\phi_i: \mathcal{V} \rightarrow \mathcal{V}$ are morphisms, and on every fiber it is true that $\sum_i(\phi_i)_t^*E_t \sim \alpha E_t$. Hence there is a function $f \in \bar{K}(\mathcal{V})^* \otimes \mathbb{R}$ and a fibral divisor $F \in $ Div$(\mathcal{V}) \otimes \mathbb{R}$ such that
\begin{center}
$\sum_i \phi_i^*E = \alpha E +$ div$(f) +F$, where $F:= \pi^*D$ for some $D$.
\end{center}
Now standard properties of local heights, for example [27] Theorem 5.4, transforms the divisorial relation above into the height relation
\begin{center}
$\sum_i \lambda_{\mathcal{V},E}(\phi_i(x),v)= \alpha \lambda_{\mathcal{V}, E}(x,v) + v(f(x)) + \lambda_{U,D} (\pi(x),v) + O(\lambda_{\partial U}(\pi(x), v)).$
\end{center}
Now we can repeat the same idea of the proof of Theorem 4.2.1 of [16], letting $\gamma(x,v) :=\sum_i \lambda_{\mathcal{V},E}(\phi_i (x),v))- \alpha \lambda_{\mathcal{V}, E}(x,v) -v(f(x))-\lambda_{U,D} (\pi(x),v) - O(\lambda_{\partial U}(\pi(x), v))$, and proceeding in the same way. This yields
\begin{center}
$\hat{\lambda}_{\mathcal{V}_t, E_t, (\mathcal{Q})_t}= \lambda_{\mathcal{V}, E}+ O(\lambda_{U,D} \circ \pi) + O(\lambda_{\partial U} \circ \pi),$
\end{center} which is almost what we want to prove. To conclude, we remember the fact that $ \sum_i(\phi_i)_t^*E_t \sim \alpha E_t$ on every fiber, so we can repeat the above argument with functions $f_1,...,f_n$ and divisors $D_1,..., D_n$ having the property that $\cap|D_i| = \emptyset.$ Then
\begin{center}
$\min \{ \lambda_{U, D_i}\}= \lambda_{U, \cap D_i}= \lambda_{U, \emptyset}$
\end{center} is $M_K$-bounded, so 
\begin{center}
$\hat{\lambda}_{\mathcal{V}_t, E_t, (\mathcal{Q})_t}= \lambda_{\mathcal{V}, E}+\min_i O(\lambda_{U,D_i} \circ \pi) + O(\lambda_{\partial U} \circ \pi)= \lambda_{\mathcal{V}, E}+ O(\lambda_{\partial U} \circ \pi).$
\end{center}
\end{proof}
{\bf Corollary 2.3:} \textit{Theorem 2.1 is true over global fields in any characteristic.}
\begin{proof}
As we have said, we just must add the previous result over all places of $K$ and use theorem 4.3.1 of [16].
\end{proof}


\section{Variation of the canonical height along sections}

In this section we have a more precise result for a one-parameter algebraic family of points. We keep almost the same notation from the previous section with the following addition:

\begin{itemize}
\item $T/K :$ we assume that the base variety $T$ has dimension 1, so $T$ is a smooth projective  curve.

\item $h_T :$ we assume that the Weil height function on $T$ corresponds to a divisor of degree 1.

\item $P:$ a section $P:T \rightarrow \mathcal{V}.$ We can think of the generic fiber $V$ of $\mathcal{V}$ as a variety over the function field $\bar{K}(T)$, and then the section $P$ corresponds to a point $P_V \in V(\bar{K}(T)).$

 \item The function field $\bar{K}(T)$ is itself an usual global height field , namely, for each point $t \in T$, there is an absolute value $\mbox{ord}_t$ on $\bar{K}(T)$ such that 

$\mbox{ord}_t(f):=$ order of vanishing of $f$ at $t$.
\end{itemize}

Further, the rational map $\phi_i: \mathcal{V} \dashrightarrow \mathcal{V}$ induces a morphism on the generic fiber $(\phi_i)_V: V \rightarrow V$, and we have $\sum_i(\phi_i)_V^* \eta_V = \alpha \eta_V,$ where $\eta_V$ is the restriction of $\eta$ to the generic fiber. \newpage This, by [16], allows us to construct the canonical height 
\begin{center}
$\hat{h}_{V, \eta_V, (\mathcal{Q})_V} : V(\overline{K(T)}) \rightarrow \mathbb{R},$ 
\end{center}
which can be evaluated at the point $P_V$. We also make $P_t=P(t)$. There are then three heights 
$\hat{h}_{V, \eta_V, (\mathcal{Q})_V}, \hat{h}_{\mathcal{V}_t, \eta_t, (\mathcal{Q})_t}$ and $h_T$ which may be compared, as Silverman did for abelian varieties in [30]. Moreover, the following theorem generalizes theorem 4.1 of [9] for the Kawaguchi canonical heights.\newline

{\bf Theorem 3.1:}
\textit{With notation as above,}
\begin{center}
$\lim_{h_T(t) \rightarrow \infty, t \in T^0(\bar{K})} \dfrac{\hat{h}_{\mathcal{V}_t, \eta_t, (\mathcal{Q})_t}(P_t)}{h_T(t)} =\hat{h}_{V, \eta_V, (\mathcal{Q})_V}(P_V).$
\end{center}

\begin{proof}
We start by stating together some results. First of all, from theorem 2.1 we have 
\begin{center}
$|\hat{h}_{\mathcal{V}_t, \eta_t, (\mathcal{Q})_t}(x)- h_{\mathcal{V}, \eta}(x)| \leq c_1h_T(t)+ c_2$ for all $t \in T^0$ and all $x \in \mathcal{V}_t.$
\end{center}
In particular, this is true for $x=P_t$, and the constants $c_1, c_2$ are independent of both $t$ and $x$. Second, we apply functoriality of Weil heights to the morphism $P:T \rightarrow \mathcal{V}.$ We note that $P$ will be a morphism, because we have assumed that $T$ is a smooth curve, so any rational map from $T$ to a variety is automatically a morphism. This gives
\begin{center}
$|h_{\mathcal{V},\eta}(P_t) - h_{T,P^*\eta}(t)| \leq c_3(P)$ for all $t \in T.$
\end{center}where $c_3(P)$ depends on the section $P$, but is independent of $t$. Third, we use [20], Chapter 3, Proposition 3.2 to describe the Weil height $h_{V, \eta_V}$ on the generic fiber in terms of intersection theory ,
\begin{center}
$|h_{V,\eta_V}(S_V) - \deg S^*\eta| \leq c_4$ for all sections $S:T \rightarrow \mathcal{V}.$
\end{center}
Fourth, we know that a canonical height is a Weil height up to a constant, and then
\begin{center}
$|\hat{h}_{V,\eta_V,(\mathcal{Q})_V}(Q_V) - h_{V,\eta_V}(Q_V)| \leq c_5$ for all $Q_V \in V(\bar{K(T)}).$
\end{center}
Using these four estimates and the triangle inequality, we compute
\begin{center}
$|\hat{h}_{\mathcal{V}_t, \eta_t, (\mathcal{Q})_t}(P_t)-\hat{h}_{V, \eta_V, (\mathcal{Q})_V}(P_V)h_T(t)| ~~~~~~~~~~~~~~~~~~~~~~~~~~~~~~~~~~~~\:\:~~~~~~~~~~~~~~\:~~~~\:\:\:\:~\:\:\:\:\:\:~~~~~~\:\:\:\:\:\:\:\:\:\:\:\:\:\:\:\:\:\:\:\:\:\:\:\:\:\:\:\:\:\:\:\:\:\:\:\:\:\:\:\:\:\:\:\:\:\:\:\:\:\:\:\:\:\:\:~~~~~\newline \newline
\leq |\hat{h}_{\mathcal{V}_t, \eta_t, (\mathcal{Q})_t}(P_t)- h_{\mathcal{V}, \eta}(P_t)|+ |h_{\mathcal{V},\eta}(P_t) - h_{T,P^*\eta}(t)| ~~~~~~~~~~~~~~~~~~~~~\:\:\:\:\:\:\:\:\:\:\:\:\:\:\:\:\:\:\:\:\:\:\:\:\:\:\:\:\:\:\:\:\:\:\:\:\:\:\:\:\:\:~~~~~~~~~~~~~~~~~~~~\newline \newline
 + |h_{T,P^* \eta }(t) -(\deg P^* \eta )h_T(t)| ~~~~~~~~~~~~~~~~~~~~~~~~~\:\:\:\:\:\:~~~~~~~~~~~~~~~~~~~~~~~~~~~~~~~~~~~~~~~~~~\:\:\:\:\:\:~~~\:\:\:~~~~~~~~~~ \:\:\:\:\:\:\:\:\:\:\:\:\:\:\:\:\:\:\:\:\:\:\:\:\:\:\:\:\:\:\:\:\:\:\:\:\:\:\:\:\:\:\:\:\:\:\:\:\:\:\:\:\:\:\:\:\:\:\:\:\:\:\:\:\:\:\:\:\:\:\newline \newline 
 + |(\deg P^* \eta)h_T(t) - h_{V, \eta_V}(P_V)h_T(t)|~~~~~~~~~~~~~~~~~~~~~~~~~~~~~~~~~~~~~~~~\:\:\:\:\:\:\:\:\:\:\:\:\:\:\:\:\:\:\:\:\:\:\:\:\:\:\:\:\:\:\:\:\:\:\:\:\:\:\:\:\:\:\:\:\:\:\:\:\:\:\:\:\:\:\:\:\:\:\:\:\:\:\:\:\:\:\:\:\:\:\:\:\:~~~~~~~~~~~~~~~~~~~~~~\:~~~~~~~~~ \newline \newline 
+ |h_{V, \eta_V}(P_V)h_T(t) - \hat{h}_{V, \eta_V, (\mathcal{Q})_V}(P_V)h_T(t)|~~~~~~~~~~~~~~~~~~~~~~~~~~~~~~~~~\:\:\:\:\:\:\:\:\:\:\:\:\:\:\:\:\:\:\:\:\:\:\:\:\:\:\:\:\:\:\:\:\:\:\:\:\:\:\:\:\:~~~~~~~~~~~~\:\:\:\:\:\:\:\:\:\:\:\:\:\:\:\:\:\:\:\:\:\:~~~~~~~~~~~~~~~~~~~ \newline \newline
\leq (c_1h_T(t) + c_2) + c_3(P) + |h_{T,P^* \eta }(t) -(\deg P^* \eta )h_T(t)| + c_4h_T(t)+ c_5h_T(t). $
\end{center}
We now divide this inequality by $h_T(t)$ and let $h_T(t) \rightarrow \infty.$ This gives
\begin{center}
$\lim \sup_{h_T(t) \rightarrow \infty}| \dfrac{\hat{h}_{\mathcal{V}_t, \eta_t, (\mathcal{Q})_t}(P_t)}{h_T(t)} - \hat{h}_{V, \eta_V, (\mathcal{Q})_V}(P_V)| \newline \newline
 \leq c_1 + c_4 + c_5 + \lim \sup_{h_T(t) \rightarrow \infty}| \dfrac{h_{T, P^* \eta}(t)}{h_T(t)}- (\deg P^* \eta)|.~~~$
\end{center}
Term $c_3(P)$ has disappeared because it depends on $P$. Moreover, Corollary 3.5 of Chapter 4 from [20] implies that the heights $h_{T, P^* \eta}$ and $(\deg P^* \eta)h_T(t)$ are quasi-equivalent, and so \begin{center}
$ \lim_{h_T(t) \rightarrow \infty}  \dfrac{h_{T, P^* \eta}(t)}{h_T(t)}= (\deg P^* \eta).$
\end{center}
This gives the fundamental estimate
\begin{center}
$\lim \sup_{h_T(t) \rightarrow \infty}| \dfrac{\hat{h}_{\mathcal{V}_t, \eta_t, (\mathcal{Q})_t}(P_t)}{h_T(t)} - \hat{h}_{V, \eta_V, (\mathcal{Q})_V}(P_V)| 
 \leq c_1 + c_4 + c_5,$
\end{center} where the constants $c_1, c_4$ and $c_5$ are independent of both the section and the point $t$, so the inequality above works with $f \circ P$ in place of $P$ for all $f \in \mathcal{Q}_n, n \in \mathbb{N}$. By (ii) of Theorem 1.2.1 from [16], we know that
\begin{center}
$\sum_{f \in (\mathcal{Q}_n)_t}\hat{h}_{\mathcal{V}_t, \eta_t, (\mathcal{Q})_t}(f(x))=\alpha^n\hat{h}_{\mathcal{V}_t, \eta_t, (\mathcal{Q})_t}(x)~~$  $\sum_{f \in (\mathcal{Q}_n)_V}\hat{h}_{V, \eta_V, (\mathcal{Q})_V}(f(x))=\alpha^n \hat{h}_{V, \eta_V, (\mathcal{Q})_V}(x).$
\end{center} 
So we finally obtain
\begin{center}
$\alpha^n\lim \sup_{h_T(t) \rightarrow \infty} |\dfrac{\hat{h}_{\mathcal{V}_t, \eta_t, (\mathcal{Q})_t}(P_t)}{h_T(t)} - \hat{h}_{V, \eta_V, (\mathcal{Q})_V}(P_V)|= ~~~~~~~~\:\:\:\:\:\:\:\:\:\:\:\:\:\:\:\:\:\:\:\:\:\:\:\:\:\:\:\:\:\:\:\:\:\:\:\:\:\:\:\:~~~~~~~~~~~~~~~~~~~~~~~~~~~~~~ \newline \newline \lim \sup_{h_T(t) \rightarrow \infty} |\dfrac{\sum_{f \in (\mathcal{Q}_n)_t}\hat{h}_{\mathcal{V}_t, \eta_t, (\mathcal{Q})_t}(f(P_t))}{h_T(t)} - \sum_{f \in (\mathcal{Q}_n)_V}\hat{h}_{V, \eta_V, (\mathcal{Q})_V}(f(P_V))| \newline \newline \leq  k^n (c_1 + c_4 + c_5).~~~~~~~~~~~~~~~~~~~~~~~~~~~~~~~~~~~\:\:\:\:\:\:\:\:\:\:\:\:\:\:\:\:\:\:\:\:\:\:\:\:\:\:\:\:\:\:\:\:\:\:\:\:\:\:\:\:\:\:\:\:\:\:\:\:\:\:\:\:\:\:\:\:\:\:\:\:\:\:\:\:\:\:\:\:\:\:\:\:\:\:\:\:\:\:\:\:\:\:\:\:\:\:\:\:\:\:\:\:\:\:\:\:\:\:\:\:\:\:\:\:\:\:\:\:\:\:\:\:\:\:~~~~~~~~~~~~~~~~~~~~~~~~~~~~~~~~~~~~~~~~~~~~~~~~~~~~~~~~~~~~~~~$
\end{center}
The right hand-side of the above inequality does not depend on $n$, while $\alpha >k$, so letting $n \rightarrow \infty$ gives us the inequality that we wanted to show.
\begin{center}
$\lim \sup_{h_T(t) \rightarrow \infty}| \dfrac{\hat{h}_{\mathcal{V}_t, \eta_t, (\mathcal{Q})_t}(P_t)}{h_T(t)} - \hat{h}_{V, \eta_V, (\mathcal{Q})_V}(P_V)| =0.$
\end{center}
\end{proof}

\section{Canonical local heights as intersection multiplicities}
In this section we show that Kawaguchi's canonical local height $\hat{\lambda}_{V,E,\mathcal{Q},}$ can be computed as an intersection number.
We fix an absolute value $v$ on $K$ and let $O_v$ denote the ring of $v$-integers in $K$. We continue with the notation used in the previous sections but we add the assumption that $V$ is a smooth projective variety, and  that the morphisms $\phi_i:V \rightarrow V$ are finite, correspondent to a dynamical system $(V, \phi_1,...,\phi_k)=(V,\mathcal{Q})$. We assume that $E$ is defined over $K$, where $E \in $ Div$(V) \otimes \mathbb{R}$ is a divisor satisfying $\sum_{i=1}^k  \phi_i^* E \sim \alpha E, \alpha > k.$
Let $S:=$ Spec$(O_v)$. We will say that a smooth scheme $\mathcal{V} /S$ is a weak N\'{e}ron Model for $(V/K, \mathcal{Q})$ over $S$ if it satisfies the following axioms: \newpage

(1) The generic fiber of  $\mathcal{V} /S$, denoted by $\mathcal{V}_K$, is $V$.

(2) Every point $P \in V(K)$ extends to a section $\bf{P}:S \rightarrow \mathcal{V}.$

(3) There exist finite morphisms $\Phi_i: \mathcal{V} /S \rightarrow \mathcal{V} /S$ whose restriction to the generic fiber are the $\phi_i.$

We note that the N\'{e}ron Model of an Abelian Variety is a weak N\'{e}ron Model for $(A/K, [n])$ for all $n \geq 2$. Indeed, for an abelian variety $A$, N\'{e}ron first showed that any canonical local height $\hat{\lambda}_{A,D}(.,v)$ can be interpreted as an intersection multiplicity on the special fiber of the N\'{e}ron model of $A$ over Spec$(\mathcal{O}_v)$(see [20], chapter 11, section 5).

Henceforth we will assume that $(V/K, \mathcal{Q})$ has a weak N\'{e}ron Model $\mathcal{V}/S.$ Let $\mathcal{V}_s$ denote the special fiber of $\mathcal{V}$ and write
\begin{center}
$\mathcal{V}_s=\sum_{j=1}^n \mathcal{V}_s^j,$
\end{center} where $\mathcal{V}_s^1,..., \mathcal{V}_s^n \in $ Div$(\mathcal{V})$ are the irreducible components of $\mathcal{V}_s.$ If $W$ is a prime divisor of $V$ rational over $K$, then $\overline{W}$, its closure in $\mathcal{V}$, is a prime divisor on $\mathcal{V}.$ Extending this process by linearity, we obtain a natural injection
\begin{center}
Div$(V)_K \rightarrow $ Div$(\mathcal{V}),~~~~~ D \rightarrow \bar{D}$ 
\end{center}
Similarly, given a point $P \in V(K)$, we write $\bar{P} = \bf{P}(S)$ to denote the image of the section $\bf{P} \in \mathcal{V}.$ Note that the divisor group on $S$ is a cyclic group generated by the special point $(s).$ Hence, for any $D \in $ Div$(V)_K$ and any $ P \in V(K)$ which does not lie in the support of $D$, we may define the intersection multiplicity $i(D,P)$ (also denoted by $\bar{P} .\bar{D}$) by 
\begin{center}
$\bf{P}^* \bar{D} =$ $ i(D,P)(s).$
\end{center}
With these notations in hand, we can now state the main result of this section, which is a more general version for results inside 6.1 of [9] due to Call and Silverman. For the proof, we will make use of a more refined  theorem in algebra linear, that was not required for the proof of their mentioned earlier result. \newline

{\bf Theorem 4.1:}
  \textit{Suppose $\mathcal{V}/S$ is a weak Neron model for $(V/K,\mathcal{Q})$ over $O_v.$ Let $\hat{\lambda}_{V,E,\mathcal{Q},f}$ be a canonical local height as constructed in theorem 4.2.1 of [16]. Moreover, suppose that $\alpha > nk$ for $n$ the number of irreducible components of the special fiber. Then there exist real numbers $\gamma_1,..., \gamma_n$ so that for all $P \in V(K) - |E|,$}
\begin{center}
$\hat{\lambda}_{V,E,\mathcal{Q},f}(P)= \bar{P} . (\bar{E} + \sum_{j=1}^n \gamma_j\mathcal{V}_s^j).$
\end{center}

An important point in the proof of this theorem is to describe the action of $\Phi_i$ on the set of irreducible components $\{ \mathcal{V}_s^1,..., \mathcal{V}_s^n \}$ of $\mathcal{V}_s$. Since $\Phi_i $ is a finite morphism, it maps each irreducible component of $\mathcal{V}_s$ onto another irreducible component (possibly the same component) of $\mathcal{V}_s$. Let $N=\{1,...,n \}.$ Then
 \begin{center}
$A_i=A_{\Phi_i}: N \rightarrow N$ defined by $\Phi_i: \mathcal{V}_s^j \rightarrow \mathcal{V}_s^{A_i(j)}$ for $j \in N$.
\end{center}
We can identify $A_i$ with a matrix of the following type. \newline

{\bf Definition 4.2:}
\textit{A square matrix $M$ is a permutation-type matrix if every column of $M$ has exactly one $1$ and all other entries are $0$.}\newline

It is a fact (lemma 6.2(b) of [9]) that every eigenvalue of a permutation matrix is 0 or is a root of unity. Such information is used in the proof of theorem 6.1 of [9]. For our more general situation, we will need the following theorem.\newline

{\bf Theorem 4.3:} \textit{Let $A$ be an $n-$square nonnegative matrix. Then}
\begin{center}
$\min_{1 \leq i \leq n} \sum_{j=1}^n a_{ij} \leq \rho(A) \leq \max_{1 \leq i \leq n} \sum_{j=1}^n a_{ij}.$
\end{center}
\textit{In other words, the spectral radius of a nonnegative square matrix is between the smallest row sum and the largest row sum.}
\begin{proof} See theorem 5.24 of [36].
\end{proof}

\begin{proof} (of Theorem 4.1). Since $E$ is assumed to be rational over $K$, we may fix a rational function $f \in K(V)^* \otimes \mathbb{R}$ so that 
\begin{center}
$ \sum_{i=1}^k \phi_i^*E= \alpha E + \mbox{div}_V(f).~~~~ (1)$
\end{center}
Since $K(V) \cong K(\mathcal{V})$, we may also regard $f$ as an element of $K(\mathcal{V})^* \otimes \mathbb{R}.$ Then the divisors of $f$ on $V$ and $\mathcal{V}$ differ by a divisor supported on the special fiber, say
\begin{center}
$\overline{\mbox{div}_V(f)}= \mbox{div}_{\mathcal{V}}(f) + Z_f,$ where $Z_f= \sum_{j=1}^n m(j,f)\mathcal{V}_s^j,~~~~ (2)$
\end{center} for some constants $m(j,f) \in \mathbb{R}.$

By Theorem 4.2.1 of [16], there is a unique canonical local height \newline $\hat{\lambda}_E:= \hat{\lambda}_{V,E,\mathcal{Q},f}$ which satisfies
\begin{center}
$\sum_{i=1}^k\hat{\lambda}_E(\phi_i (p),v)= \alpha\hat{\lambda}_E(p,v) + v(f(p)). ~~~~(3)$
\end{center}

Consider the map $V(K) - |E| \rightarrow \mathbb{R}$ defined by $P \rightarrow i(E,P) = \bar{P} . \bar{E}.$ Given any $P \in V(K) -|E|,$ there is a pair $(U,g)$ representing $\bar{E}$ such that $U \subset \mathcal{V}$ is an open neighborhood of $P$ and $g(P) \neq 0, \infty.$ Then, by definition, $i(E,P)= v(g(P)),$ independent of the choice of the pair $(U,g).$ Thus, the map $P \mapsto i(E,P)$ is a Weil Local Height function for $E$ on $V(K)$.

Note that $\Phi_i^* \bar{E}$ and $\overline{\phi_i^*E}$ differ by a divisor supported on the special fiber, since $\Phi_i$ and $\phi_i$ are the same on the generic fiber. Combining this fact with (1), we have 
\begin{center}
$\sum_i \Phi_i^* \bar{E} = \alpha \bar{E} + \overline{\mbox{div}_V(f)} + \sum_{j=1}^n n_j \mathcal{V}_s^j,~~~~~~~(4)$
\end{center} for some constants $n_j \in \mathbb{R}$. Further, \begin{center} $(\Phi_i)_* \bar{P}=(\Phi_i)_* \bf{P}(S)=$$\Phi_i $$\circ \bf{P}(S)= \bf{\phi}_i(P)$$(S) =\overline{\phi_i(P)}$, \end{center} where $\bf{\phi_i(P)}$ is the section corresponding to $\phi_i (P).$ Hence, 
\begin{center}
$ \sum_i \bar{P} . \Phi_i^* \bar{E} = \sum_i (\Phi_i)_* \bar{P} . \bar{E} = \sum_i \overline{\phi_i(P)} . \bar{E} = \sum_i i(E,\phi_i(P)). ~~~~~~~~~ (5)$
\end{center}
Intersecting both sides of (2) with $\bar{P}$ yields: \newline

$\bar{P} . \overline{\mbox{div}_V(f)}= \bar{P} . \mbox{div}_{\mathcal{V}}(f) + \bar{P} . Z_f = v(f(P)) + \bar{P} . \sum_{j=1}^n m(j,f) \mathcal{V}_s^j. ~~~~(6)$ \newline

Now, intersecting both sides of (4) with $\bar{P}$ and applying (5) and (6), we conclude
\begin{center}
$\sum_i i(E, \phi_i(P))= \alpha i(E,P) + v(f(P)) + \bar{P} . \sum_{j=1}^n c_j \mathcal{V}_s^j, ~~~~~~~~~~~ (7)$
\end{center}
where $c_j= m(j,f) + n_j$ are constants which depend on $E, \mathcal{Q}$ and $f$, but are independent of $P.$ In particular, we see that (7) holds for all $P \in V(K)$ for which the intersection multiplicities $i(E, \phi_i(P))$ and $i(E,P)$ are defined, i.e, for all $P \notin |E| \cup |\phi_1^*E| \cup ... \cup |\phi_k^* E|.$

Next, we will show that one can choose real numbers $x_1,..., x_n$ so that the function 
\begin{center}
$\Lambda_E(P) = i(E,P) + \bar{P} . \sum_{j=1}^n x_j \mathcal{V}_s^j ~~~~~~ (8)$
\end{center} satisfies
\begin{center}
$\sum_i \Lambda_E(\phi(P))= \alpha \Lambda_E(P) + v(f(P)) ~~~~~~~~~(9)$

\end{center}
For all $P \in V(K) -(|E| \cup |\phi_1^*E| \cup ... \cup |\phi_k^* E|).$ Using (8) and (7), we compute

$\newline \sum_i \Lambda_E ( \phi_i(P)) - \alpha \Lambda_E(P) - v(f(P))
\newline \newline  ~~~~~~~~= \sum_{i=1}^k (\overline{\phi_i(P)} . \sum_{j=1}^n x_j \mathcal{V}_s^j) - \alpha \bar{P} .  \sum_{j=1}^n x_j \mathcal{V}_s^j + \bar{P} .  \sum_{j=1}^n c_j \mathcal{V}_s^j.$ \newline \newline
Recall that $\Phi_i$ determines a permutation type matrix $A_i:= A_{\Phi_i}$ defined by $\Phi_i: \mathcal{V}_s^j \rightarrow \mathcal{V}_s^{A_i(j)}$. Since $\bar{P}$ and $\overline{\phi_i(P)}= \Phi_i(\bar{P})$ intersect the components of $\mathcal{V}_s$ transversally, it follows from the definition of $A_i$ that if $\bf{P}(s) \in \mathcal{V}_s^t$, then
\begin{center}
$\bar{P} . \sum_{j=1}^n x_j \mathcal{V}_s^j= x_t$ and $ \overline{\phi_i(P)}.  \sum_{j=1}^n x_j \mathcal{V}_s^j= x_{A_i(t)}.~~~~~$  (10)
\end{center}
Therefore it suffices to find constants $x_1,..., x_n$ such that
\begin{center}
$\sum_{i=1}^k x_{A_i(t)} - \alpha x_t + c_t = 0$ for $t=1,...,n.$
\end{center}
Writing $x_1,...,x_n$ and $c_1,..., c_n$ as column forms, we can combine these $n$ equations into a matrix equation
\begin{center}
$(\alpha \bf{I} - \sum_{i=1}^kA_i) \bf{x} = \bf{c} . $
\end{center} The $A_i$ are permutation-type matrices, so Theorem 3.7 says that the absolute value of an eigenvalue of $\sum_i A_i$ is at most $nk$, but $\alpha >nk$ by hypotheses. So $\det(\alpha \bf{I} - \sum_{i=1}^kA_i) \neq 0 $ and $(\alpha \bf{I} - \sum_{i=1}^kA_i)$ is invertible and we may take $\bf{x} = (\alpha \bf{I} - \sum_{i=1}^kA_i)^{-1} \bf{c}$. This finishes the proof that we can choose $x_1,...,x_n$ so that the function $\Lambda$ defined by (8) satisfies (9).

To complete the proof, we will show that $\hat{\lambda}_E(P,v) = \Lambda_E(P)$ for all $P$ in $V(K) - |E|.$ Since $\hat{\lambda}_E(.,v)$ and $i(E,.)$ are both Weil local heights for $E$, their difference has a unique $v$-continuous extension to a bounded $v$-continuous function defined on all of $V(K)$ (see [20], chapter 10, proposition 1.5, and 2.3). Hence, by (8), we see that the map $L_E(P):= \hat{\lambda}_E(P,v)- \Lambda_E(P)$ extends to a bounded function on $V(K)$, namely, by a constant C $\geq 0$. Further, since $\hat{\lambda}_E$ and $\Lambda_E$ satisfy (3) and (9), it follows that
\begin{center}
$\sum_i L_E(\phi_i(P))= \alpha L_E(P)$ for all $ P \in V(K)$.
\end{center}
Therefore, for any $P \in V(K),$
\begin{center}
$|L_E(P)| \leq |\alpha^{-N} \sum_{\phi \in \mathcal{Q}_N} L_E(\phi(P))| \leq \frac{k^N}{\alpha^N}. C \rightarrow_{N \rightarrow \infty} 0.$
\end{center}
We conclude that $L_E\equiv0$, so $\hat{\lambda}_E(P,v)= \Lambda_E(P)$ $\forall P \in V(K) - |E|.$
\end{proof}

\section{Canonical Metrics of Commuting Systems}

With two polarized dynamical systems $(X,\mathcal{F}=\{f_1,...,f_k\}, \mathcal{L}, \alpha), \alpha >k$ and $(X,\mathcal{G}=\{g_1,...,g_t\}, \mathcal{L}, \beta), \beta>t$, we can build two canonical metrics, two canonical heights, and two canonical measures for $\mathcal{L} \in ~$Pic$(X) \otimes \mathbb{R}$. We will see that if all the maps of one of the systems commute with all the maps of the other, then the canonical metrics, canonical heights, and canonical measures associated to each system are identical.

We consider a projective variety $X$ over a number field $K$, and $(X;f_1,...,f_k)$ a dynamical eigensystem of $k$ morphisms $\mathcal{F}:=\{f_1,...,f_k\}$ over $K$ associated with an ample line bundle $\mathcal{L} \in $ Pic$(X) \otimes \mathbb{R}$ of degree $\alpha>k$ as in section 3 of [16], so we have an isomorphim $\phi: \mathcal{L}^{\otimes^{\alpha}} \rightarrow^{\cong} f_1^{*}\mathcal{L} \otimes ... \otimes f_k^{*}\mathcal{L} $. This situation will be called a polarized dynamical eigensystem $(X, f_1,..., f_k, \mathcal{L}, \alpha)$ on $X$ defined over $K$.
Assume that for every place $v$ of $K$ we have chosen a continuous and bounded metric $||.||_v$ on each fibre of $\mathcal{L}_v:= \mathcal{L} \otimes_K K_v$. The following theorem is stated in theorem 3.3.1 of [16] for $||.||_{\infty}$. \newline

{\bf Theorem 5.1:} \textit{The sequence defined recurrently by $||.||_{v,1}:=||.||_v$ and \begin{center}$||.||_{v,n}=(\phi^{*}(f_1^{*}||.||_{v,n-1}...f_k^{*}||.||_{v,n-1}))^{\frac{1}{\alpha}}$ for $n>1$ \end{center}
converge uniformly on $X(\bar{K}_v)$ to a metric $||.||_{v,\mathcal{F}}$ (independent of the choice of $||.||_{v,1}$) on $\mathcal{L}_v$ which satisfies the equation }
\begin{center} $||.||_{v,\mathcal{F}}=(\phi^{*}(f_1^{*}||.||_{v,\mathcal{F}}...f_k^{*}||.||_{v,\mathcal{F}}))^{\frac{1}{\alpha}}$.
\end{center}
\begin{proof} The proof is the same as is theorem 3.3.1 of [16] with $v$ in place of $\infty$.
\end{proof}

{\bf Definition 5.2:} \textit{The metric $||.||_{v,\mathcal{F}}$ is called the canonical metric on $\mathcal{L}_v$ relative to the system of maps $\mathcal{F}.$} \newline

The following proposition relates the canonical metrics associating to commuting maps. It represents the main result of this section, and it is a natural and simple generalization of proposition 2.5 of [26] for the metric defined in [16, theorem 3.1.1].\newline

{\bf Proposition 5.3:}
\textit{Let $(X,\mathcal{F}=\{f_1,...,f_k\}, \mathcal{L}, \alpha)$ and $(X,\mathcal{G}=\{g_1,...,g_t\}, \mathcal{L}, \beta)$ be two polarized systems with $\alpha >k, \beta >t$ on $X$ defined over $K$. Suppose that $f_i \circ g_j= g_j \circ f_i $ for all $i,j$. Then $||.||_{v,\mathcal{F}}=||.||_{v,\mathcal{G}}$.}
\begin{proof} The key idea is that the canonical metric does not depend on the metric from which we have started the iteration with. Let $s \in \Gamma (X, \mathcal{L})$ be a non-zero section of $\mathcal{L}$. We are going to consider two metrics $||.||_{v,1}=||.||_{v,\mathcal{F}}$ and $||.||^{`}_{v,1}=||.||_{v,\mathcal{G}}$ on the line bundle $\mathcal{L}$. By our definition of canonical metric for $\mathcal{F}$, we can start with 
$||.||^{`}_{v,1}$ and obtain $||s(x)||_{v,\mathcal{F}}= \lim_r (||\prod_{f \in \mathcal{F}_r}s(f(x))||^ `_{v,1})^{\frac{1}{\alpha^{r}}}$, but also by our definition of canonical metric for $\mathcal{G}$ starting with \begin{center} $||.||_{v,1}=||.||_{v,\mathcal{F}}$ we get $||s(x)||_{v,\mathcal{G}}= \lim_l ||\prod_{g \in \mathcal{G}_l}s(g(x))||^{\frac{1}{\beta^{l}}}_{v,\mathcal{F}}$.\end{center} So using the uniform convergence and the commutativity of the maps,
\begin{center}
$||s(x)||_{v,\mathcal{F}}=\lim_r ||\prod_{f \in \mathcal{F}_r}s(f(x))||^{\frac{1}{\alpha^{r}}}_{v,\mathcal{G}} ~~~~~~~~~~~~~~~~~~~~~~~~~~~\:\:\:~~\newline \newline ~~~~~= \lim_{r,l} \prod_{f \in \mathcal{F}_r}\prod_{g \in \mathcal{G}_l}||s(f(g(x)))||^{\frac{1}{\alpha^{r} \beta^{l}}}_{v,1} ~~~~~~~~~~~~~~~~~~~~~~~~~~~~~~~~~~~~~\newline \newline ~~~~~
=  \lim_{r,l} \prod_{f \in \mathcal{F}_r}\prod_{g \in \mathcal{G}_l}||s(g(f(x)))||^{\frac{1}{\alpha^{r} \beta^{l}}}_{v,1}~~~~~~~~~~~~~~~~~~~~~~~~~~~~~~~~~~~~~ \newline \newline=\lim_l ||\prod_{g \in \mathcal{G}_l}s(g(x))||^{\frac{1}{\beta^{l}}}_{v,\mathcal{F}}
~~~~~~~~~~~~~~~~~~~~~~~~~~~~~~~~~~~\:\:\:\:\:\:\:\:\:\:\:\:\:\:\:\:\:\:\:\:\:~\:~~\:~~~~~\newline \newline
= ||s(x)||_{v,\mathcal{G}}$,~~~~~~~~~~~~~~~~~~~~~~~~~~~~~~~~~~~~~\:\:\:\:\:\:\:\:\:\:\:\:\:\:\:\:\:\:\:\:\:\:\:\:\:\:\:\:\:\:\:\:\:\:\:\:\:\:\:\:\:\:\:\:\:\:\:\:\:\:\:\:\:\:\:\:\:\:\:\:\:\:\:\:\:\:\:\:\:\:\:\:\:\:\:\:\:\:\:~~~~~~~~~~~~~~~~~~~~~~~~~~~~~~~~~~~~~~
\end{center} which was the result we wanted to prove.
\end{proof}

Let $X$ $n$-dimensional projective variety defined over a number field $K$ and $(X,\mathcal{F}=\{f_1,...,f_k\}, \mathcal{L}, \alpha)$ a polarized system with $\alpha >k$ defined over $K$. Let $v$ be a place of $K$ over infinity. We can consider morphisms \begin{center}$f_{i} \otimes v:X_v \rightarrow X_v$ over the complex variety $X_v$.\end{center} Associated to $\mathcal{F}$ and $v$ we also have the canonical metric $||.||_{v,\mathcal{F}}$ and therefore the distribution $\hat{c}_1(\mathcal{L},||.||_{v,\mathcal{F}})= \frac{1}{\pi i} \partial\bar{\partial}\log ||s||_{v,\mathcal{F}}$, with $  s \in \Gamma (X, \mathcal{L})- \{0\}$, analogous to the first Chern form of $(\mathcal{L},||.||_{v,\mathcal{F}})$. It can be proved that this is a positive current in the sense of Lelong, as in section 3.2 of [16], we can define the product
\begin{center}
$\hat{c}_1(\mathcal{L},||.||_{v,\mathcal{F}})^{n}=\hat{c}_1(\mathcal{L},||.||_{v,\mathcal{F}})...\hat{c}_1(\mathcal{L},||.||_{v,\mathcal{F}})$,
\end{center} which represents a measure $\mu$ on $X_v$.\newline

{\bf Definition 5.4:} \textit{The measure $d\mu_{\mathcal{F}}=\hat{c}_1(\mathcal{L},||.||_{v,\mathcal{F}})^{n}/\mu(X)$, is called the canonical measure associated to $\mathcal{F}$ and $v$. Once we fix $\mathcal{L}$, it depends only on the metric $||.||_{v,\mathcal{F}}$.}

Now we assume that $X$ is an arithmetic variety of absolute dimension $n+1$, that is, given a number field $K$, $X$ is flat and of finite type over Spec$\mathcal({O}_K)$ of relative dimension $n$. We can define (see section 2 of [16]) the arithmetic intersection number $\hat{c}_1(\mathcal{L}_1)...\hat{c}_1(\mathcal{L}_{n+1})$ of the classes of the hermitian line bundles $(\mathcal{L}_i,||.||)$ on $X$, which means that each line bundle $\mathcal{L}_i$ on $X$ is equipped with a hermitian metric $||.||_{v,i}$ over $X_v= X \otimes_K$ Spec$(\mathcal{O}_K)$, for each place $v$ at infinity. Such line bundles are called adelic line metrized bundles when they can be equipped with semipositive metric for all places $v$. So we can define adelic intersection numbers $\hat{c}_1({\mathcal{L}_1}_{|Y})...\hat{c}_1({\mathcal{L}_{n+1}}_{|Y})$ over a $p$-cycle $Y \subset X$. Suppose that we are in the presence of a polarized dynamical eigensystem $(X,\mathcal{F}, \mathcal{L}, \alpha)$. In this case the canonical metric $||.||_{\mathcal{F}}$ represents a semipositive metric on $\mathcal{L}$, and we can define the canonical height associated to $(\mathcal{L}, ||.||_{\mathcal{F}})$.\newline

{\bf Definition 5.5:}
\textit{The canonical height $\hat{h}_{\mathcal{F}}(Y)$ of a $p$-cycle $Y$ in $X$ is defined as}
\begin{center}
$\hat{h}_{\mathcal{F}}(Y)= \dfrac{\hat{c}_1(\mathcal{L}_{|Y})^{p+1}}{(\dim Y +1)c_1(\mathcal{L}_{|Y})^{p}}$.
\end{center} It depends only on $(\mathcal{L},||.||_{\mathcal{F}})$, where $||.||_{\mathcal{F}}$ is actually representing a colection of canonical metrics over all places of $K$. An inportant particular case of canonical height will be the canonical height $\hat{h}_{\mathcal{F}}(P)$ of a point $P$ of $X$. Since the canonical measure and the canonical height were defined depending only on the canonical metric of the system, the equality of canonical metrics and measure is a direct consequence of proposition 5.3, as we just state below.\newline

{\bf Proposition 5.6:} \textit{Let $(X,\mathcal{F}=\{f_1,...,f_k\}, \mathcal{L}, \alpha)$ and $(X,\mathcal{G}=\{g_1,...,g_t\}, \mathcal{L}, \beta)$ be two polarized systems with $\alpha >k, \beta >l$ on $X$ defined over $K$. Suppose that $f_i \circ g_j= g_j \circ f_i $ for all $i,j$. Then $\hat{h}_{\mathcal{F}}=\hat{h}_{\mathcal{G}}$ and $d\mu_{\mathcal{F}}=d\mu_{\mathcal{G}}$.}
\begin{proof}
It is a consequence of the last two definitions that the statements depend only on the canonical metric presented in proposition 5.3.
\end{proof}


\begin{thebibliography}{9}                                                                                                %


\bibitem{} F. Amoroso and S. David \textit{ Distribution des points de petite hauteur dans les groupes multiplicatifs}. Ann. Scuola Norm. Cl. Sci. (5) 3, no. 2, 325-348

\bibitem{} F. Amoroso and S. David, \ \textit{Minoration de la hauteur normalise\'{e}e dans un tore}, J. Inst. Math. Jussieu 2 (2003), no. 3, 335-381.

\bibitem{}  F. Amoroso and E. Viada, \ \textit{Small points on subvarieties of a torus}, Duke Math. J. 150 (2009), no. 3, 407-442 .

\bibitem{} A. Baker, \ \textit{New Advances in Transcendence Theory}, Cambridge University Press(1988).

\bibitem{}  A. Baragar,\ \textit{Canonical vector heights on algebraic K3 surfaces with Picard number two}, Canad. Math. Bull. 46 (2003), 495-508 .

\bibitem{}  A. Baragar,\ \textit{Rational points on K3 surfaces in $\mathbb{P}^1 \times\mathbb{P}^1 \times \mathbb{P}^1$ },Math. Ann. 305 (1996),541-558 .

\bibitem{}  Y. Bilu,\ \textit{ Limit distribution of small points on algebraic tori}, Duke Math. J. 89 (1997), 465-508.

\bibitem{}  E. Bombieri and W. Gubler,\ \textit{ Heights in Diophantine Geometry}, Number 4 in New Mathematical Monographs. Cambridge University Press, Cambridge, 2006.

\bibitem{} G. Call and H. Silverman, \ \textit{Canonical heights on varieties with morphisms},
Compositio Math. 89 (1993), 163-205. 

\bibitem{}  A. Chambert-Loir,\ \textit{ Mesures et \'{e}quidistribution sur les espaces de Berkovich}, J. Reine Angew. Math., 595 (2006), 215-325 .

\bibitem{}  W. Fulton, \ \textit{Intersection theory, 2nded.,Ergeb.Math.Grenzgeb.(3)2}, Springer-Verlag, Berlin 1998.

\bibitem{}  V.Guedj, \ \textit{Ergodic properties of rational mappings with large topological degree}, Ann.of Math.(2) 161 (2005), no. 3, 1589-1607.

\bibitem{}  R. Hartshorne, \ \textit{Algebraic Geometry}, Springer-Verlag, New
York, 1977.

\bibitem{}  M. Hindry and J. Silverman, \ \textit{Diophantine Geometry: An Introduction}, volume 201 of Graduate Texts in Mathematics. Springer-Verlag, New
York, .

\bibitem{}  S. Kawaguchi, \ \textit{Canonical height functions for affine plane automorphisms}, Math. Ann., 335(2):285-310, 2006.

\bibitem{}  S. Kawaguchi, \textit{Canonical heights, invariant currents, and dynamical eigensystems of morphisms for line bundles,} preprint, math.NT/0405006. (The revised version is accepted in J. Reine Angew. Math.)

\bibitem{}  S. Kawaguchi, \ \textit{Projective surface automorphisms of positive topological entropy from an arithmetic viewpoint},  available at arXiv:math/0510634v2 [math.AG].

\bibitem{}  S. Kawaguchi and Joseph H. Silverman, \ \textit{On the dynamical and arithmetic
degrees of rational self-maps of algebraic varieties},J. Reine Angew. Math. 713 (2016),
21–48.

\bibitem{} S. Kawaguchi and J. H. Silverman, \ \textit{Dynamical canonical heights for Jordan blocks, arithmetic degrees of orbits, and nef canonical heights on abelian varieties}, preprint 2013, http://arxiv.org/abs/1301.4964

\bibitem{}  S. Lang, \ \textit{Fundamentals of Diophantine Geometry},  New York, (1983).

\bibitem{}  C. G. Lee, \ \textit{Equidistribution of periodic points of some automorphisms on K3 surfaces},  available at arXiv:1102.4860v1[math.NT].

\bibitem{} J. Lin and C. Wang, \ \textit{Canonical height functions for monomial maps},2012. arXiv:1205.2020v1.

\bibitem{}  A. Moriwaki,\ \textit{ Arithmetic height functions over finitely generated fields}, Invent. Math. 140 (2000),101-142 .

\bibitem{}  A. Moriwaki,\ \textit{The canonical arithmetic height of subvarieties of an abelian variety over a finitely generated field }, J. reine angrew. Math. 530 (2001), 33-54.

\bibitem{}  B. and M. Vaughan , \ \textit{Jordan Normal and Rational Normal Form Algorithms},2004. arXiv:cs/0412005v1.

\bibitem{}  J. Pineiro, \ \textit{Canonical metrics of commuting maps}, available at arXiv:0803.3468v1[math.NT].

\bibitem{}  J. Silverman, \ \textit{Arithmetic distance functions and height function in Diophantine
 geometry},  Math. Ann. 279(1987) 193-216.

\bibitem{}  J. Silverman, \ \textit{Examples of dynamical degree equals arithmetic degree}, available at arXiv:1212.3015v1[math.NT].

\bibitem{}  J. Silverman, \ \textit{Dynamical Degrees, Arithmetic Entropy, and Canonical Heights for Dominant Rational Self-Maps of Projective Space}, available at arXiv:1111.5664[math.NT].

\bibitem{}  J. Silverman, \ \textit{Heights and the specialization maps for families of abelian varieties},  J. Reine Angew. Math. 342(1983) 197-211.

\bibitem{}  J. Silverman, \ \textit{Rational points on K3 surfaces: a new canonical height}, Invent. Math. 105 (1991), 347-373.

\bibitem{}  J. Silverman, \ \textit{The Arithmetic of Dynamical Systems}, volume 241 of Graduate Texts in Mathematics. Springer, New York, 2007.

\bibitem{} J. Silverman, \ \textit{The arithmetic of elliptic curves}, 2nd ed., Grad. Texts in Math. 106, Springer-Verlag, Dordrecht 2009.

\bibitem{} M. Waldschmidt, \ \textit{Diophantine approximation
on linear algebraic groups. Transcendence properties of the exponential function
in several variables}, Springer-Verlag, Berlin, 2000.

\bibitem{}  L. Wang, \ \textit{Rational points and canonical heights on K3-surfaces in $\mathbb{P}^1 \times\mathbb{P}^1 \times \mathbb{P}^1$ }, Recent developments in the inverse Galois problem (Seattle, WA, 1993), Contemp. Math. 186 (1995), 273-289.

\bibitem{}  F. Zhang, \ \textit{Matrix Theory: Basic Results and Techniques}, Springer-Verlag, (2011).

\bibitem{}  S. Zhang, \ \textit{Small points and adelic metrics}, J. Alg. Geom. 4 (1995),281-300 .





\end{thebibliography}
\end{document}